\documentclass[10pt, amsmath, amssymb]{article} 
\usepackage{epsfig, amssymb}
\newcommand{\finpr}{\hfill $\square$ \vspace{2mm}}

\newcommand{\C}{{\mathbb{C}}}

\def\be{\begin{eqnarray}}
\def\ee{\end{eqnarray}}
\def\bee{\begin{eqnarray*}}
\def\eee{\end{eqnarray*}}
\newtheorem{thm}{Theorem}
\newtheorem{cor}{Corollary}
\newtheorem{lem}{Lemma}

\newtheorem{prop}{Proposition}

\begin{document}

\title{Edge-local equivalence of graphs}
\author{Maarten Van den Nest\footnote{Corresponding author. E-mail: mvandenn@esat.kuleuven.be}, \  Bart De Moor,\\ \ \\
ESAT-SCD, K.U. Leuven,\\  Kasteelpark Arenberg 10,\\
 B-3001 Leuven, Belgium}

\date{\today}
\def\makeheadbox{}

\maketitle

\begin{abstract}
The local complement $G*i$ of a simple graph $G$ at one of its
vertices $i$ is obtained by complementing the subgraph induced by
the neighborhood of $i$ and leaving the rest of the graph
unchanged. If $e=\{i,j\}$ is an edge of $G$ then $G*e=((G*i)*j)*i$
is called the edge-local complement of $G$ along the edge $e$. We
call two graphs edge-locally equivalent if they are related by a
sequence of edge-local complementations. The main result of this
paper is an algebraic description of edge-local equivalence of
graphs in terms of linear fractional transformations of adjacency
matrices. Applications of this result include (i) a polynomial
algorithm to recognize whether two graphs are edge-locally
equivalent, (ii) a formula to count the number of graphs in a
class of edge-local equivalence, and (iii) a result concerning the
coefficients of the interlace polynomial, where we show that these
coefficients are all even for a class of graphs; this class
contains, as a subset, all strongly regular graphs with parameters
$(n, k, a, c)$, where $k$ is odd and $a$ and $c$ are even.
\end{abstract}


\section{Introduction}

Let $G$ be a simple graph with vertex set $V$ and edge set $E$.
The local complement of $G$ at a vertex $i\in V$, denoted by
$G*i,$ is obtained by replacing the subgraph induced by the
neighborhood of $i$ by its complement, and leaving the rest of the
graph unchanged. Two graphs are called locally equivalent if they
are related by a sequence of local complementations. If $e=\{i,
j\}\in E$, the edge-local complement $G*e$ of $G$ along the edge
$e$ is defined by $G*e=((G*i)*j)*i$. We call two graphs
edge-locally equivalent if they are related by a sequence of
edge-local complementations.

Local equivalence of graphs was studied in detail in the 1980s,
primarily by Bouchet \cite{Bouchet, alg_bouchet, Bou_di,
Bou_circle, Bou_trees}, where, among other results, this work lead
to a polynomial algorithm to recognize circle graphs. Bouchet
studied local equivalence of graphs in the framework of certain
algebraic and combinatorial structures called \emph{isotropic
systems,} which can be regarded as subgroups of
$K\times\dots\times K$ that are self-dual with respect to some
bilinear form, where $K$ is the Klein four group \cite{Bou_iso,
Bou_graph_iso, Bou_conn_iso}. Recently, local equivalence of
graphs has found a renewed interest in a number of distinct fields
of research. First, local equivalence of graphs turns out to be a
relevant equivalence relation in the context of quantum
information theory\footnote{We note that quantum information
theory is the field of research of the present authors.}, where
graphs correspond to an important class of (pure) states of
distributed quantum systems, called \emph{graph states}
\cite{entgraphstate}; in this context local equivalence of graphs
corresponds to so-called \emph{local Clifford equivalence}, or LC
equivalence, of such graph states, an equivalence relation that is
of importance in the study of multi-partite entanglement in graph
states \cite{localcliffgraph, LU_LC, cliff_inv_prl},
measurement-based quantum computation \cite{1wayQC} and the theory
of quantum error-correction \cite{Gott, Glynn}.  Second, there is
an intimate connection between (edge-)local complementation of
graphs and the \emph{interlace} polynomial \cite{ABS00, BBCP02,
BBRS01, AH04, interlace_prop}, which is a recently introduced
graph polynomial motivated by questions arising from DNA
sequencing by hybridization \cite{ABC00}; the interlace polynomial
can be defined recursively through identities involving edge-local
complementation. Third, LC equivalence of graph states, and thus
local equivalence of graph, has recently been studied in classical
cryptography in the theory of (quadratic) generalized bent
functions \cite{bent_I, bent_II}.

Whereas local complementation of graphs is understood reasonably
well, a general theory for edge-local equivalence seems to be
missing. It is the aim of this paper to take the first significant
steps towards this end. In the following we perform an algebraic
analysis of edge-local complementation of graphs, relating this
notion to certain transformations over GF(2) of adjacency
matrices, namely \emph{linear fractional transformations} (LFTs)
having the form \be\Gamma\mapsto ((A+I)\Gamma +
A)(A\Gamma+A+I)^{-1},\ee where $\Gamma$ is the adjacency matrix of
a graph on $n$ vertices and $A$ is an $n\times n$ diagonal matrix.
Our main result is a proof that two graphs are edge-locally
equivalent if and only if their adjacency matrices are related by
an LFT of the above type. We will then use this result in a number
of applications: first, a polynomial time algorithm to recognize
edge-local equivalence of two arbitrary graphs is obtained; the
complexity of this algorithm is ${\cal O}(n^4)$, where $n$ is the
number of vertices of the considered graphs. Second, the
description of edge-local equivalence in terms of LFTs allows to
obtain a formula to count the number of graphs in a class of
edge-local equivalence. Third, we consider a number of invariants
of edge-local equivalence. Fourth, we prove a result concerning
the coefficients of the interlace polynomial of a graph, where we
show that these coefficient are all even for a class of graphs;
this class contains, among others, all strongly regular graphs
with parameters $(n, k, a, c)$, where $k$ is odd and $a$ and $c$
are even. Fifth,  we show that a nonsingular adjacency matrix can
be inverted by performing a sequence of edge-local
complementations on the corresponding graph. Finally, we consider
the connection between edge-local equivalence of graphs and the
action of local Hadamard matrices on the corresponding graph
states, yielding a restricted version of LC equivalence of graph
states. We believe that this last result is of interest in the
theory of generalized bent functions.

This paper is organized as follows. in section \ref{sectie_2} we
introduce the basic definitions regarding local complementation
and edge-local complementation. In section \ref{sectie_3} we
present the basic theory involving linear fractional
transformations of adjacency matrices and we state our main
result, which is subsequently proven in section \ref{sectie_4}.
Then in section \ref{sectie_5} a number of applications of this
result are considered.

\subsection*{Notations}

We note that in this paper \emph{all} arithmetic is considered
over the field $\mathbb{F}_2=$ GF(2), except in section 5.5, where
also complex numbers are involved.

The following notations will be used.  The identity matrix is
denoted by $I$, and $E_i$ is a square matrix where $(E_i)_{ii}=1$
and all other entries are zero (the dimensions of $I$ and $E_i$
follow from the context). Let $n\in\mathbb{N}_0$. If
$\omega\subseteq\{1, \dots, n\}$, then $\bar\omega$ denotes the
complement of $\omega$ in $\{1, \dots, n\}$. Let $X$ be an
$n\times n$ matrix. For every $\omega\subseteq\{1, \dots, n\}$,
define the $|\omega|\times |\omega|$ principal submatrix
$X[\omega]$ of $X$ by\be X[\omega]=(X_{ij})_{i,j\in\omega},\ee and
define the $|\omega|\times (n-|\omega|)$ off-diagonal submatrix
$X\langle\omega\rangle$ by  \be
X\langle\omega\rangle=(X_{ij})_{i\in\omega, j\in\bar\omega}.\ee
Also, for every $n\times n$ diagonal matrix $A$, we denote $\hat
A:=(A_{11}, \dots, A_{nn})\in\mathbb{F}_2^n$. For every $v\in
\mathbb{F}_2^n$, $\hat v$ is the $n\times n$ diagonal matrix with
the entries of $v$ on the diagonal. Finally, we denote the support
of the vector $v$ and its associated diagonal matrix by
\be\mbox{supp}(v) = \mbox{ supp}(\hat v) = \{i\in\{1, \dots, n\}\
|\ v_i=1\}.\ee

In this paper all graphs are finite, undirected and simple (no
loops, no multiple edges). A graph $G$ with vertex set $V$ and
edge set $E$ is denoted by $G=(V, E)$ and, for simplicity and
without loss generality, we will always take $V=\{1, \dots, n\}$
for some $n\in \mathbb{N}_0$. The adjacency matrix
$\Gamma(G)\equiv \Gamma$ of the graph $G$ is the unique $n\times
n$ $0-1$ matrix such  that $\Gamma_{ij}=1$ if and only if $\{i,
j\}\in E$. Adjacency matrices are symmetric and
have zeros on their diagonals. 

The neighborhood $N(i) \subseteq V$ of a vertex $i\in V$ is the
set of all vertices that are adjacent to $i$. For a subset $\omega
\subseteq V$ of vertices, the induced subgraph $G[\omega]
\subseteq G$ is the graph with vertex set $\omega$ and edge set
\be\{\{i, j\} \in E\ |\ i, j \in \omega\}.\ee Note that the
adjacency matrix of $G[\omega]$ is $\Gamma[\omega]$. The
complement $\bar G$ of $G=(V, E)$ is the graph on the same vertex
set $V$ with the property that $\{i, j\}$ is an edge of $\bar G$
if and only if it is not an edge of $G$.

\section{Edge-local equivalence}\label{sectie_2}

\begin{figure}
\hspace{1cm}\includegraphics[width=10cm]{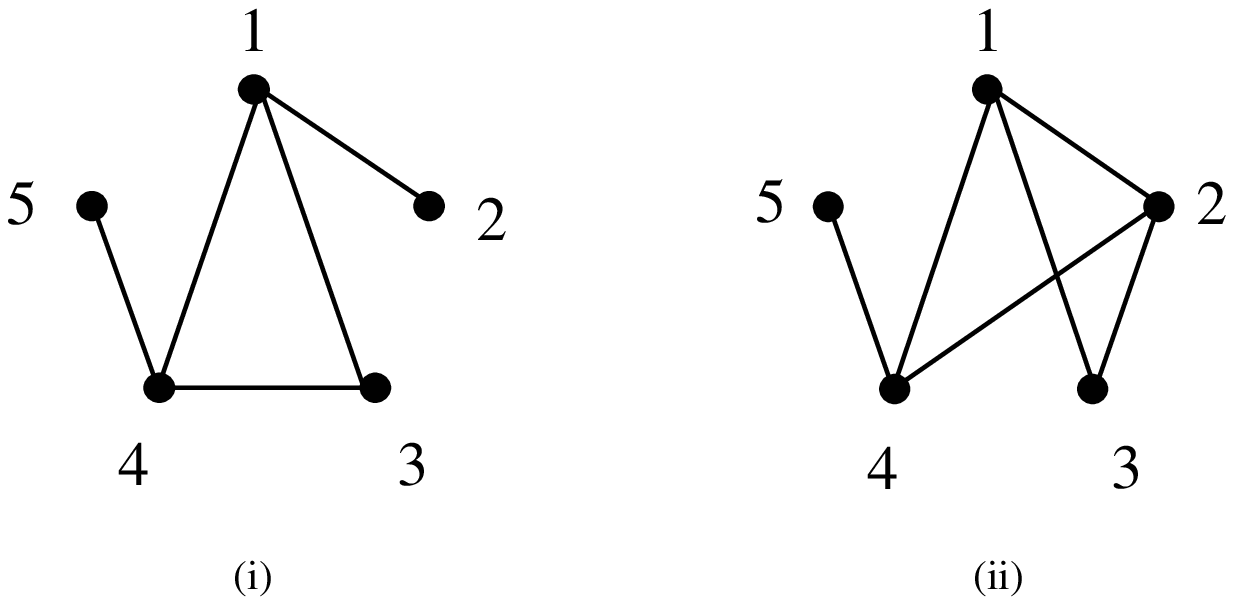}
\caption{\label{fig:LC} (i) Graph on 5 vertices, and (ii) its
local complement at vertex $1$. The neighborhood of vertex 1 in
(i) consists of the vertices 2, 3 and 4. Hence, the graph in (ii)
is obtained by complementing the induced subgraph on these
vertices, and leaving the rest of the graph unchanged.}
\end{figure}

Let $G=(V, E)$ be a graph with adjacency matrix $\Gamma$ and let
$i\in V$ be an arbitrary vertex of $G$. The \emph{local complement
of $G$ at vertex $i$}, denoted by $G*i$, is the graph on the same
vertex set $V$ obtained by replacing the subgraph $G[N(i)]$ by its
complement and leaving the rest of the graph unchanged. Note that
$(G*i)*i=G$. One can easily verify that the adjacency matrix
$\Gamma*i$ of $G*i$ is equal to \be \Gamma*i = \Gamma +
\Gamma_i\Gamma_i^T + \hat{\Gamma_i},\ee where $\Gamma_i$ is the
$i$th column of $\Gamma$. An example of local complementation is
given in Fig. \ref{fig:LC}.

The local complementation rule gives rise to an equivalence
relation on the set of graphs, called \emph{local equivalence}.
Two graphs $G$ and $G'$ on the same vertex set $V$ are
 called \emph{locally equivalent} if there exist $i_1,\dots, i_N\in V$
 such that \be G' =(((G*i_1)*i_2)\dots)*i_N. \ee

Local complementation is also used as an elementary building block
for a second graph transformation, which we call \emph{edge-local
complementation}. Letting $G=(V, E)$ be a graph and $e=\{i, j\}\in
E$, one defines the \emph{edge-local complement of $G$ along the
edge $e$}, denoted by $G*e$, by \be G*e = ((G*i)*j)*i =
((G*j)*i)*j.\ee More explicitly, the graph $G*e$ is the unique
graph with the following properties.

\begin{itemize}

\item Let $k, l\in V\setminus e$. Then $k$ is adjacent to $i$
($j$) in $G*e$ if and only if it is adjacent $j$ ($i$) in $G$.

\item Let $k, l\in V\setminus e$ and suppose that one of the
following situations occurs:
    \begin{itemize}
    \item[(i)] $k\in N(i)\setminus N(j)$ and $l\in N(j)\setminus N(i)$;
    \item[(ii)] $k\in N(j)\setminus N(i)$ and $l\in N(i)\setminus N(j)$;
    \item[(iii)] $k\in N(i)\setminus N(j)$ and $l\in N(i)\cup N(j)$;
    \item[(iv)] $k\in N(j)\setminus N(i)$ and $l\in N(i)\cup N(j)$;
    \item[(v)] $k\in N(i)\cup N(j)$ and $l\in N(i)\setminus N(j)$;
    \item[(vi)] $k\in N(i)\cup N(j)$ and $l\in N(j)\setminus N(i)$;
    \end{itemize}
    then $\{k, l\}$ is an edge of $G*e$ if and only if $\{k, l\}$ is
    not an edge of $G$. If none of the cases (i)-(vi) occur, then
    $\{k, l\}$ is an edge of $G*e$ if and only if it is
    an edge of $G$.
\end{itemize}
\begin{figure}
\hspace{4cm}\includegraphics[width=4cm]{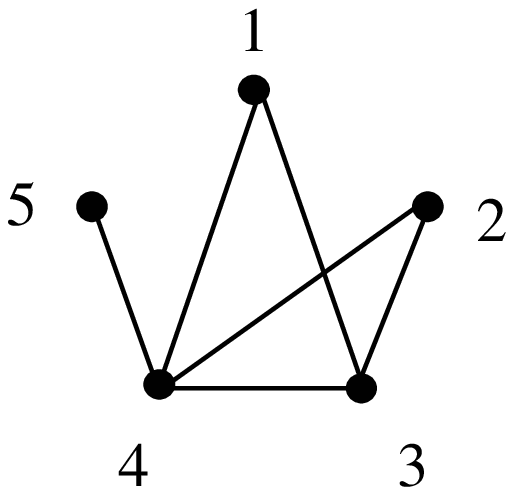}
\caption{\label{fig:ELC} Edge-local complement of graph (i) in
Figure 1 along the edge $\{1, 3\}$, which can be obtained by first
complementing at vertex 1, then at vertex 3, and again at vertex
1.}
\end{figure}
Note that $(G*e)*e=G$. An example of local complementation along
an edge is given in Fig. \ref{fig:ELC}. Note that local
complementation of $G$ along an edge $e$ can also be formulated
elegantly in terms of adjacency matrices. Letting $\Gamma$ and
$\Gamma*e$ be the adjacency matrices of $G$ and $G*e$,
respectively, the matrix $\Gamma*e$ is obtained as follows:
 \be\label{LCE_alg} \begin{array}{lcl}(\Gamma*e)[e] &=& \left[ \begin{array}{cc} 0&1\\1&0 \end{array}\right]\\
(\Gamma*e)\langle e\rangle &=&
\left[ \begin{array}{cc} 0&1\\1&0 \end{array}\right]\Gamma\langle e\rangle\\
(\Gamma*e)[\bar e] &=& \Gamma[\bar e] + \Gamma \langle
e\rangle^T\left[
\begin{array}{cc} 0&1\\1&0
\end{array}\right]\Gamma\langle e\rangle.\end{array} \ee
We will need this representation in section 4.

Analogous to the case where local complementation is used to
define local equivalence of graphs, edge-local complementation
gives rise to an equivalence relation as well, which we call
\emph{edge-local equivalence}.  Letting $[V]^2$ be the set of all
2-element subsets of $V$, two graphs $G$ and $G$ are called
\emph{edge-locally equivalent} if there exist $e_1, \dots, e_N\in
[V]^2$ such that
\begin{itemize}
\item[(i)] $e_1$ is an edge of $G$ and $e_{\alpha}$ is an edge of
$((G*e_1)*\dots)*e_{\alpha-1}$, for every $\alpha=2, \dots, N$,
and \item[(ii)]  $ ((G*e_1)*\dots)*e_N = G'$.
\end{itemize}

\section{Linear fractional transformations}\label{sectie_3}

In this section we consider matrix transformations of the form
\be\Gamma\mapsto (A\Gamma + B)(C\Gamma + D)^{-1},\ee where are
$\Gamma$, $A$, $B$, $C$ and $D$ are $n\times n$ matrices. These
transformations are called linear fractional transformations
(LFTs). In particular, we will consider the case where $\Gamma$ is
the adjacency matrix of a graph and $A$, $B$, $C$ and $D$ are
diagonal matrices satisfying $AD+BC=I$; this last constraint is
equivalent to stating that the $2n\times 2n$ matrix
\be\label{Q_def} Q:=\left[
\begin{array}{cc}A&B\\C&D\end{array}\right]\ee is nonsingular. The
set of all such matrices $Q$ is a group isomorphic to $GL(2,
\mathbb{F}_2)^n$ and will be denoted by ${\cal C}_n$. We will also
write $Q=[A, B, C, D]$ instead of (\ref{Q_def}). If $\Gamma$ is an
$n\times n$ adjacency matrix, $Q=[A, B, C, D]\in{\cal C}_n$ and
$C\Gamma+D$ is nonsingular, we will use the notation \be Q(\Gamma)
= (A\Gamma + B)(C\Gamma + D)^{-1}.\ee Note that $Q(\Gamma)$ is
symmetric whenever $\Gamma$ is symmetric. However, it is possible
that $Q(\Gamma)$ does not have zero diagonal when $\Gamma$ does.
Therefore, we associate to every $Q=[A, B, C, D]\in{\cal C}_n$ a
domain of definition $\Delta(Q)$, consisting of all $n\times n$
adjacency matrices $\Gamma$ such that
\begin{itemize}
\item[(i)] $C\Gamma+D$ is nonsingular, and \item[(ii)] $Q(\Gamma)$
has zero diagonal.
\end{itemize}

There is an intimate connection between local complementation of
graphs and the above LFTs of adjacency matrices, as we showed in
Ref. \cite{localcliffgraph}, where we studied local
complementation in the context of quantum information theory. The
main result in \cite{localcliffgraph} is the following:

\begin{thm} Let $G$, $G'$ be graphs with $n\times n$ adjacency matrices
$\Gamma$ and $\Gamma'$, respectively. Then $G$ and $G'$ are
locally equivalent if and only if there exists a matrix $Q\in
{\cal C}_n$ such that $\Gamma\in\Delta(Q)$ and
$Q(\Gamma)=\Gamma'$.
\end{thm}
For example, letting $\Gamma$ be an $n\times n$ adjacency matrix
and $i\in \{1, \dots, n\}$, and defining $Q^{\Gamma}_i = [I,
\hat{\Gamma_i}, E_i, I]$, where $\Gamma_i$ is the $i$th column of
$\Gamma$, one has \be \Gamma*i = Q_i^{\Gamma}(\Gamma).\ee We note
that the result in theorem 1 is also implicit in the work of
Bouchet regarding isotropic systems and local complementation
\cite{Bouchet}. We elaborate on this matter in appendix A. Glynn
also proved an equivalent of theorem 1 independently of our own
work \cite{Glynn}.

It is the aim of this paper to derive a result analogous to
theorem 1 for edge-local complementation of graphs. To do so, we
will consider LFTs having the following specific structure. Let
${\cal H}_n$ be the subset of ${\cal C}_n$ consisting of all
elements of the form \be H=[A+I, A, A, A+I]=:H^A\ee for some
diagonal $n\times n$ matrix $A$. Note that ${\cal H}_n$ is a
subgroup of ${\cal C}_n$. Also, $H^AH^B = H^{A+B}$ for all
diagonal matrices $A$ and $B$, such that ${\cal H}_n$ is
isomorphic to the additive group of $\mathbb{F}_2^n$. The main
result of this paper can now be formulated as follows:
\begin{thm}
Let $G$ and $G'$ be graphs with adjacency matrices $\Gamma$ and
$\Gamma'$, respectively. Then $G$ and $G'$ are edge-locally
equivalent if and only if there exists an operator $H\in {\cal
H}_n$ such that $\Gamma\in\Delta(H)$ and $H(\Gamma)=\Gamma'$.
\end{thm}
The proof of this result will be given in section 4. For the
remainder of this section, we analyze some basic properties of
LFTs that we will need below. A first basic result is the
following:
\begin{prop}
Let $\Gamma$, $\Gamma'$ be $n\times n$ adjacency matrices and let
$Q=[A, B, C, D]\in{\cal C}_n$. Then the following statements are
equivalent:
\begin{itemize}
\item[(i)] $\Gamma\in\Delta(Q)$ and $Q(\Gamma)=\Gamma'$;
\item[(ii)]$\label{lastbasic}\Gamma'C\Gamma + A\Gamma + \Gamma'D +
B=0.$
\end{itemize}
\end{prop}
{\it Proof:} The implication from (i) to (ii) is trivial. We prove
the reverse implication. First, define the $n-$dimensional linear
space \be\label{V_Gamma} V_{\Gamma} = \{(\Gamma u,u)\ |\
u\in\mathbb{F}_2^n\}\subseteq\mathbb{F}_2^{2n}\ee and the space
$V_{\Gamma'}$ analogously. Note that the columns of the matrix
$[\Gamma\ I]^T$ ($[\Gamma'\ I]^T$) are a basis of $V_{\Gamma}$
($V_{\Gamma'}$). Second, consider the following (symplectic) inner
product $\langle\cdot, \cdot\rangle$ on $\mathbb{F}_2^{2n}$: \be
\langle (u,v), (u', v')\rangle = u^Tv' + v^Tu',\ee where $u, v,
u', v'\in\mathbb{F}_2^n$. It is easy to verify that $\langle x,
x\rangle = 0$ for every $x\in V_{\Gamma}$, such that this space is
self-dual with respect to the above inner product, and the same
holds for $V_{\Gamma'}$. Furthermore, note that the identity (ii)
holds if and only if \be
\left[\begin{array}{cc} I &\Gamma' \end{array}\right] Q \left[\begin{array}{c}\Gamma\\
I\end{array}\right] =0,\ee which is in turn equivalent to stating
that $\langle x, y\rangle = 0$ for every $x\in V_{\Gamma'}$ and
$y\in QV_{\Gamma}$. Since the space $V_{\Gamma'}$ is self-dual, it
follows that $QV_{\Gamma} = V_{\Gamma'}$. Therefore, the columns
of the matrix $Q[\Gamma\ I]^T$ and those of the matrix $[\Gamma'\
I]^T$ form bases of the same linear space, such that there must
exist a nonsingular matrix $R$ satisfying $Q[\Gamma\ I]^T =
[\Gamma'\ I]^TR$. More explicitly, this last equation reads \be
\left[
\begin{array}{c} A\Gamma + B\\ C\Gamma+D\end{array}\right] = \left[
\begin{array}{c} \Gamma'R\\ R\end{array}\right],\ee showing that
$\Gamma\in\Delta(Q)$ and $Q(\Gamma) = \Gamma'$. This completes the
proof.\finpr

\noindent The above result will be used in section 5.1 to obtain
an efficient algorithm to recognize edge-local equivalence of two
given graphs.

Next we state a criterion to verify whether an adjacency matrix
belongs to the domain of a given element of ${\cal C}_n$.

\begin{prop}\label{prop_3_3}
Let $Q=[A, B, C, D]\in{\cal C}_n$ with $\omega:= \mbox{ supp}(C)$
and let $\Gamma$ be an $n\times n$ adjacency matrix. Then
$\Gamma\in \Delta(Q)$ if and only if $\Gamma[\omega] +D[\omega]$
is nonsingular and $\Gamma \hat{AC}=\hat{BD}$.
\end{prop}
{\it Proof:} First we show that $C\Gamma +D$ is invertible if and
only if $\Gamma[\omega] +D[\omega]$ is invertible. It follows from
$Q\in {\cal C}_n$ that $D[\bar\omega]=I$. Therefore, $C\Gamma +D$
has the form \be\left[
\begin{array}{c|c} I& 0\\\hline * & \Gamma[\omega] + D[\omega]\end{array}\right]\ee
up to a simultaneous permutation of its rows and columns. It is
now easy to see that $C\Gamma +D$ is nonsingular if and only
$\Gamma[\omega] +D[\omega]$ is, which proves the first part of the
proposition.

Second, we show that the matrix $(A\Gamma + B)(C\Gamma + D)^{-1}$
has zeros on its diagonal if and only if $\Gamma
\hat{AC}=\hat{BD}$. To see this, note that for every $n\times n$
adjacency matrix $X$ and $R\in GL(n, \mathbb{F}_2)$ the matrix
$R^T X R$ is symmetric and has zero diagonal. This follows from
the observation that $(R^T X R)_{ii}= R_i^TX R_i=0$, where $R_i$
is the ith column of $R$ (we have $R_i^TX R_i=0$ since $X$ is
symmetric and has zero diagonal). It follows from the above that
$(A\Gamma + B)(C\Gamma + D)^{-1}$ has zero diagonal if and only if
$\left(\Gamma C +D \right)\left( A\Gamma + B\right)$ has, which is
equivalent to \be\label{diagzero}\left(\Gamma AC\Gamma + AD\Gamma
+ \Gamma BC + BD \right)_{ii}=0\ee for every $i=1, \dots, n$.
Since $\Gamma$ has zero diagonal one has $(AD\Gamma)_{ii}= (\Gamma
BC)_{ii}=0$. Moreover, \be(\Gamma AC\Gamma)_{ii}&=&\sum_{k=1}^n
\Gamma_{ik}A_kC_k \Gamma_{ki}\nonumber\\
&=& \sum_{k=1}^n \Gamma_{ik}^2A_kC_k \nonumber\\&=& \sum_{k=1}^n
\Gamma_{ik}A_kC_k= \Gamma_i^T \hat{AC},\ee where $\Gamma_i$ is the
$i$th column of $\Gamma$. This shows that (\ref{diagzero}) is
equivalent to $\Gamma \hat{AC}=\hat{BD}$. This proves the result.
\finpr

\begin{prop}
Let $\Gamma$ be an $n\times n$ adjacency matrix and let $Q=[A, B,
C, D], Q'=[A', B', C', D']\in{\cal C}_n$ such that
$\Gamma\in\Delta(Q)$ and $Q(\Gamma)\in \Delta(Q')$. Then
$\Gamma\in \Delta(Q'Q)$ and $Q'(Q(\Gamma)) = (Q'Q)(\Gamma)$.
\end{prop}
{\it Proof:} denoting $\Gamma':= Q'(Q(\Gamma))$, it follows from
proposition 1 that $\Gamma\in \Delta(Q'Q)$ and
$(Q'Q)(\Gamma)=\Gamma'$ if and only if \be\Gamma'V\Gamma + T\Gamma
+ \Gamma'W + U = 0, \ee where \be\left[ \begin{array}{cc} T&U\\
V&W\end{array}\right] = Q'Q= \left[
\begin{array}{cc} A'A+B'C&A'B+B'D\\ C'A+D'C&C'B+D'D\end{array}\right].\ee
Substituting \be\Gamma' =
(A'(A\Gamma+B)(C\Gamma+D)^{-1}+B')(C'(A\Gamma+B)(C\Gamma+D)^{-1}+D')^{-1}\ee
yields the result after a straightforward calculation. \finpr

\begin{prop} Let $\Gamma$ be an $n\times n$ adjacency matrix and let $H^A\in{\cal H}_n$ with $\omega = \mbox{ supp}(A)$.
 Then the following statements are equivalent.
 \begin{itemize}
 \item[(i)] $\Gamma\in\Delta(H^A)$;
 \item[(ii)] $A\Gamma+A+I$ is nonsingular;
 \item[(iii)] $\Gamma[\omega]$ is nonsingular.
 \end{itemize}
\end{prop}
{\it Proof:} the implication from (i) to (ii) is trivial. The
implication from (i) to (iii) follows from proposition 2. The
reverse implication is proven by using proposition 2 and noting
that $A(A+I)=0$. The equivalence of (ii) and (iii) is proven by an
argument analogous to the (the first paragraph of the) proof of
proposition 2. This completes the proof. \finpr

\begin{prop}
Let $\Gamma$ be an $n\times n$ adjacency matrix and $H^A\in{\cal
H}_n$ with $\Gamma\in\Delta(H^A)$. Let $\omega=\mbox{ supp}(A)$.
Then $H^A(\Gamma)$ has the following
structure: \be H^A(\Gamma)[\omega] &=& \Gamma[\omega]^{-1}\\
H^A(\Gamma)\langle\omega\rangle &=&
\Gamma[\omega]^{-1}\Gamma\langle\omega\rangle\\
H^A(\Gamma)[\bar\omega] &=& \Gamma[\bar\omega] +
\Gamma\langle\omega\rangle^T\Gamma[\omega]^{-1}\Gamma\langle\omega\rangle.\ee
\end{prop}
{\it Proof:} the proof is obtained by a straightforward
calculation.\finpr

\begin{prop}
Let $\Gamma$ be an $n\times n$ adjacency matrix and let $A$ and
$B$ be $n\times n$ matrices such that $\Gamma\in\Delta(H^A)$ and
$H^A(\Gamma)\in \Delta(H^B)$. Then $\Gamma\in \Delta(H^{A+B})$ and
$H^B(H^A(\Gamma)) = H^{A+B}(\Gamma)$.
\end{prop}
{\it Proof: } This is an immediate corollary of proposition 3.
\finpr

\section{Edge-local equivalence of graphs and LFTs} \label{sectie_4}

In this section it is our aim to prove theorem 2. Letting $G$ be a
graph with $n\times n$ adjacency matrix $\Gamma$ and letting
$e=\{i, j\}$ be an edge of $G$, it immediately follows from
(\ref{LCE_alg}) and proposition 5 that
 \be \Gamma*e = H^{E_i+E_j}(\Gamma).\ee
Moreover, according to proposition 6 successively complementing
$G$ along edges can also be realized as an LFT. Indeed, let $e_1,
\dots, e_N\in [V]^2$ such that $e_1$ is an edge of $G$ and
$e_{\alpha}$ is an edge of $((G*e_1)*\dots)*e_{\alpha-1}$, for
every $\alpha=2, \dots, N$. Denoting $e_{\alpha}=\{i_{\alpha},
j_{\alpha}\}$ and \be A = \sum_{\alpha=1}^N\ E_{i_{\alpha}} +
E_{j_{\alpha}}, \ee proposition 6 shows that
$\Gamma\in\Delta(H^A)$ and \be ((\Gamma*e_1)*\dots)*e_N =
H^A(\Gamma).\ee   This proves the forward implication of theorem
2:

\begin{prop}
Let $G$ and $G'$ be edge-locally equivalent graphs with $n\times
n$ adjacency matrices $\Gamma,\ \Gamma'$, respectively. Then there
exists a matrix $H\in {\cal H}_n$ such that $\Gamma\in\Delta(H)$
and $\Gamma'=H(\Gamma)$.
\end{prop}

We now prove the reverse implication of theorem 2. To do so,
define the set ${\cal T}_n'$ by \be {\cal T}_n' = \{A\Gamma + A+I\
|\ A \mbox{ is } n\times n \mbox{ diagonal, } \Gamma\mbox{ is an }
n\times n\mbox{ adjacency matrix}\},\nonumber \ee and let ${\cal
T}_n:={\cal T}_n'\cap GL(n, \mathbb{F}_2)$. For every $i=1, \dots
n$ consider the transformation $f_i: \mathbb{F}_2^{n\times n}\to
\mathbb{F}_2^{n\times n}$ defined by \be f_i(X)= X(E_i X + X_{ii}
E_i + I),\ee for every $n\times n$ matrix $X$. We also denote
$f_{ij}:= f_if_jf_i$. We are now ready to state the following
lemma.
\begin{lem}
Let $R= A\Gamma + A+I\in {\cal T}_n$. Then there exist
$i_{\alpha}, j_{\alpha}\in \mbox{ supp}(A)$, where $\alpha=1,
\dots, N$, such that $f_{i_Nj_N}\dots f_{i_1j_1}(R)=I$  and \be
1=R_{i_1j_1} = \left(f_{i_1j_1}(R)\right)_{i_2j_2} = \dots=
(f_{i_{N-1}j_{N-1}}\left(\dots f_{i_1j_1}(R))\right) _{i_Nj_N}.\ee

\end{lem}
{\it Proof:} For every $i\in \{1, \dots, n\}$, let $e_i$ be the
$i$th canonical basis vector of $\mathbb{F}_2^n$. Fix
$i_1\in\mbox{ supp}(A)$ and apply $f_{i_1}$ to $R$. It can easily
be verified that the diagonal of $f_{i_1}(R)$ is equal to diag$(
R) + {R} _{i_1}$, where ${R} _{i_1}$ is the $i_1$th column of $R$
and diag$(R)$ is the diagonal of $R$. Since $R$ is nonsingular,
${R} _{i_1}$ has some nonzero element, say $R_{j_1i_1}=1$ (and one
moreover has $j_1\in\mbox{ supp}(A)$). Therefore, the
$j_1j_1-$entry of $f_{i_1}(R)$ is equal to 1. Now we apply
$f_{j_1}$, and one can verify that the $j_1$th row of
$f_{j_1}f_{i_1}(R)$ is equal to $e_{j_1}^T$. Furthermore, the
$i_1i_1-$entry of $f_{j_1}f_{i_1}(R)$ is $R_{i_1j_1}$, where, from
the symmetry of $R$, one has $R_{i_1j_1}=R_{j_1i_1}=1$. By again
applying $f_{i_1}$, one finds that the $i_1$th row of
$f_{i_1j_1}(R)$ is equal to $e_{i_1}^T$ and the $j_1$th row
remains equal to $e_{j_1}^T$. Moreover, one can also verify that
\be f_{i_1j_1}(R) = A'\Gamma' + A'+I\in {\cal T}_n, \ee for some
$n\times n$ adjacency matrix $\Gamma'$ and some diagonal matrix
$A'$ with $\mbox{supp}(A') \subseteq\mbox{supp}(A)\setminus \{i_1,
j_1\}$. The rest of the proof follows by induction on
$|\mbox{supp}(A)|$. This proves the result. \finpr

\noindent The above result will now be used to complete the proof
of theorem 2.

\begin{prop}
Let $\Gamma$ be an $n\times n $ adjacency matrix and let $A$ be a
diagonal matrix such that $\Gamma\in\Delta(H^A)$. Let $i_{\alpha},
j_{\alpha}\in \mbox{ supp}(A)$, for every $\alpha=1, \dots, N$, be
the sequence of vertices produced by applying lemma 1 to $A\Gamma
+ A+I$. Then, denoting $e_{\alpha}=\{i_{\alpha}, j_{\alpha}\}$,
one has
\begin{itemize} \item[(i)] $e_1$ is an edge of $G$ and
$e_{\alpha}$ is an edge of $((G*e_1)*\dots)*e_{\alpha-1}$, for
every $\alpha=2, \dots, N$, and \item[(ii)]  $
((\Gamma*e_1)*\dots)*e_N = H^A(\Gamma)$,
\end{itemize} implying that $\Gamma$ and
$H^A(\Gamma)$ are edge-locally equivalent.
\end{prop}
{\it Proof:} The result is proven by induction on $N$. We start
with the basis of the induction, i.e., $N=1$, where we consider a
diagonal matrix $A$ such that $\Gamma\in\Delta(H^A)$ and
\be\label{basis} f_{ij}(A\Gamma+A+I)=I.\ee for some $i, j\in\mbox{
supp}(A)$. Denoting $e=\{i, j\}$, the identity (\ref{basis})
implies that $R=A\Gamma+A+I$ satisfies \be R[e] = \left[
\begin{array}{cc} 0&1\\1&0 \end{array}\right],\quad R\langle \bar
e\rangle =0,\quad  R[\bar e] = I, \ee showing that $A = E_i+ E_j$.
Recalling that $H^{E_i+E_j}(\Gamma)=\Gamma*e$ proves the basis of
the induction.

In the induction step, we assume that the result is true for
sequences $f_{j_1k_1}, \dots, f_{j_Nk_N}$ of length $N$ and prove
that this implies the result is true for sequences of length
$N+1$. Thus, we consider a diagonal matrix $A$ such that
$\Gamma\in\Delta(H^A)$ and, denoting $R=A\Gamma+A+I$, we assume
that \be
 \label{induct}f_{i_Nj_N} \dots
f_{i_1j_1}f_{ij}(R)=I\ee and \be 1=R_{ij} =
\left(f_{ij}(R)\right)_{i_1j_1} = \dots=
(f_{i_{N-1}j_{N-1}}\left(\dots f_{ij}(R))\right) _{i_Nj_N},\ee for
some $i, j, i_{\alpha}, j_{\alpha}\in \mbox{ supp}(A)$ for every
$\alpha=1, \dots, N$. Denoting $e=\{i,j\}$, we then define
\begin{eqnarray}
R' &=& f_{ij}(R)\nonumber\\
G'&=& G*e \nonumber\\
 \Gamma'&=&\Gamma*e \nonumber\\
A' &=& A + E_i + E_j\end{eqnarray} Note that $R_{ij}=1$ implies
that $e$ is an edge of $G$, such that $G*e$ is well-defined. It is
now straightforward to show that $R'=A'\Gamma'+A'+I$, which is the
crucial observation of the proof. We then have the following
situation:
\begin{itemize}
\item $\Gamma'\in \Delta(H^{A'})$: this follows from the
invertibility of $A'\Gamma'+A'+I = R'$ and proposition 4, and
\item $f_{i_Nj_N} \dots f_{i_1j_1}(R')=I$
\end{itemize}The induction hypothesis then implies that
\begin{itemize} \item[(i)] $e_1$ is an edge of $G'$ and
$e_{\alpha}$ is an edge of $((G'*e_1)*\dots)*e_{\alpha-1}$, for
every $\alpha=2, \dots, N$, and \item[(ii)]  $
((\Gamma'*e_1)*\dots)*e_N = H^{A'}(\Gamma')$.
\end{itemize}
Moreover, we have $\Gamma'=\Gamma*e= H^{E_i+E_j}(\Gamma)$, such
that \be H^{A'}(\Gamma') = H^{A'}H^{E_i+E_j}(\Gamma) =
H^{A}(\Gamma).\ee This completes the proof of the proposition.
\finpr

The proof of theorem 2 is now obtained by combining propositions 8
and 9.

\section{Applications} \label{sectie_5}

In this section we present some applications of theorem 2. First,
we present an efficient algorithm to recognize whether two given
graphs are edge-locally equivalent; second, we derive a formula to
count the number of graphs in a class of edge-local equivalence;
third, we consider some graph invariants under edge-local
equivalence; fourth, we use theorem 2 to obtain some results
regarding the coefficients of the interlace polynomial of a graph;
fifth, we show that a nonsingular adjacency matrix can be inverted
by performing a sequence of edge-local complementations on the
corresponding graph; finally, we state an equivalent version of
theorem 2 in terms of graph states and local Hadamard
transformations.

\subsection{Recognizing edge-local equivalence efficiently}

Theorem 2 allows to efficiently recognize whether two given graphs
are edge-locally equivalent. To see this, let $G$ and $G'$ be two
graphs with $n\times n$ adjacency matrices $\Gamma$ and $\Gamma'$,
respectively. Then $G$ and $G'$ are edge-locally equivalent if and
only if there exists an $n\times n$ diagonal matrix $A$ such that
$H^A(\Gamma)=\Gamma'$. Following proposition 1, this happens if
and only if \be\label{system}\Gamma'A\Gamma + (A+I)\Gamma +
\Gamma'(A+I) + A=0,\ee or, equivalently,
\be\label{system'}(\Gamma'+I)A(\Gamma + I)=\Gamma+\Gamma'.\ee
Regarding the adjacency matrices $\Gamma$ and $\Gamma$ as given
and the entries of the matrix $A$ as unknowns, this is a system of
$n^2$ affine equations in $n$ unknowns, which can be solved in
$O(n^4)$ time. The graphs $G$ and $G'$ are edge-locally equivalent
if and only if (\ref{system'}) has a solution; moreover, if a
solution $A$ is found, the following algorithm produces a sequence
of edge-local complementations transforming $G$ into $G$': first,
using the method presented in the proof of lemma 1, one obtains
$i_{\alpha}, j_{\alpha}\in \mbox{ supp}(A)$, where $\alpha=1,
\dots, N$, such that \be f_{i_Nj_N}\dots
f_{i_1j_1}(A\Gamma+A+I)=I\ee Denoting $e_{\alpha} = \{i_{\alpha},
j_{\alpha}\}$, it follows from proposition 9 that \be
((\Gamma*e_1)*\dots)*e_N = H^A(\Gamma) = \Gamma',\ee yielding the
desired sequence of edge-local complementation transforming $G$
into $G'$.

Note that the system of linear equations (\ref{system'}) can be
written more explicitly as follows. Letting $A=\hat a$, where
$a=(a_1, \dots, a_n)\in\mathbb{F}_2^n$, and denoting by
$\tilde\Gamma_i$ ($\tilde\Gamma_i'$) the $i$th column of
$\Gamma+I$ ($\Gamma'+I$), for every $i=1, \dots, n$, the identity
(\ref{system'}) is equivalent to \be \label{system''}\sum_{i=1}^n\
a_i \tilde\Gamma_i'\ \tilde\Gamma_i^T=\Gamma+\Gamma'.\ee Now we
let $\gamma$ ($\gamma'$) be the $n^2-$dimensional vector which is
obtained by assembling the entries of $\Gamma$ ($\Gamma'$) row per
row in a vector. The reader can now verify that (\ref{system''})
is equivalent to \be \label{system'''}\sum_{i=1}^n\ a_i
\tilde\Gamma_i'\otimes \tilde\Gamma_i=\gamma+\gamma'.\ee Here
$\otimes$ denotes the tensor product, or Kronecker product, of
vectors. Note that (\ref{system'''}) is an equality of
$n^2-$dimensional vectors, whereas (\ref{system''}) was an
equality of $n\times n$ matrices. Thus, we have shown that $G$ and
$G'$ are edge-locally equivalent if and only if the system
\be\label{tilde_Gamma} \left[
\begin{array}{ccc}\tilde\Gamma'_1\otimes\tilde\Gamma_1&\dots&
\tilde\Gamma'_n\otimes\tilde\Gamma_n\end{array}\right] a =
\gamma+\gamma'\ee has a solution.

\subsection{Number of graphs in a class of edge-local equivalence}

A second application of theorem 2 is a formula to count the number
of graphs edge-locally equivalent to a given one. Let $G$ be a
graph with $n\times n$ adjacency matrix $\Gamma$ and let $L_e(G)$
be the set of all graphs edge-locally equivalent to $G$. Define
the following two sets \be\Delta_e(G)&:=&\{x\in \mathbb{F}_2^n\ |\
\Gamma\in\Delta(H^{\hat x})\}\nonumber\\
\Sigma_e(G)&:=&\{x\in \mathbb{F}_2^n\ |\ H^{\hat
x}(\Gamma)=\Gamma\}.\ee The set $\Delta_e(G)$ collects all
operations in ${\cal H}_n$ that have $\Gamma$ in their domain. The
set $\Sigma_e(G)$, which is a subset of $\Delta_e(G)$, collects
all operations in ${\cal H}_n$ that have $\Gamma$ as a fixed
point. Using theorem 2, it is clear that \be\label{aantal}
|L_e(G)| = \frac{|\Delta_e(G)|}{|\Sigma_e(G)|}.\ee First, it
follows from proposition 4 that \be\label{nonsing} |\Delta_e(G)|
=|\{ \omega\subseteq \{1, \dots, n\}\ |\ \Gamma[\omega] \mbox{ is
nonsingular}\}|.\ee Thus, $|\Delta_e(G)|$ simply counts the number
of nonsingular principal submatrices of $\Gamma$. Second, using
proposition 1 we find that $\Sigma_e(G)$ consists of all
$x\in\mathbb{F}_2^n$ satisfying \be\label{sigma} 0 = \Gamma
\hat{x}\Gamma + \hat{x}\Gamma + \Gamma \hat{x} + \hat{x} = (\Gamma
+ I)\hat{x}(\Gamma +I).\ee Note that $\Sigma_e(G)$ is a linear
vector space of dimension at most $n$, such that
$|\Sigma_e(G)|=2^l$, where $l$ be calculated efficiently. The
discussion at the end of section 5.1 shows that the space
$\Sigma_e(G)$ is the null space of the $n^2\times n$ matrix
\be\label{tilde_Gamma'} \left[
\begin{array}{ccc}\tilde\Gamma_1\otimes\tilde\Gamma_1&\dots& \tilde\Gamma_n\otimes\tilde\Gamma_n\end{array}\right].\ee
This implies that the dimension of $\Sigma_e(G)$ is equal to the
corank of the $n^2\times n$ matrix (\ref{tilde_Gamma'}).

Further, we recall the definition of the \emph{bineighborhood
space} of $G$, as defined in \cite{alg_bouchet}: using the
notation $\nu_{ij}:= (\Gamma_{i1}\Gamma_{j1}, \dots,
\Gamma_{in}\Gamma_{jn}) \in\mathbb{F}_2^n,$ for every $i, j\in\{1,
\dots, n\}$, the bineighborhood space $\nu(G)$ of $G$ is  the
subspace of $\mathbb{F}_2^n$ spanned by the sets of vectors \be\{
\sum_{ (i,j)\in C} \nu_{ij}\ |\ C \mbox{ is a cycle of }
G\}\mbox{\quad and \quad} \{\nu_{ij}\ |\ (i,j)\notin E \}.\ee We
can then formulate the following properties.
\begin{prop}
Let $G$ be a graph. Then \be \Sigma_e(G)\subseteq \mbox{
ker}(\Gamma +I)\cap \nu(G)^{\perp}.\ee
\end{prop}
{\it Proof:} Let $x\in\Sigma_e(G)$. Then by definition \be
&&(\Gamma_i+e_i)x=((\Gamma+I)\hat x(\Gamma+I))_{ii}=0\label{sigma_1}\\
&&\sum_{k=1}^n\Gamma_{ik}\Gamma_{jk}x_k + (x_i+x_j)\Gamma_{ij}
=0\label{sigma_2}\ee for every $i,j=1, \dots, n$, $i\neq j$, where
$e_i$ is the $i$th canonical basis vector of $\mathbb{F}_2^n$.
Equation (\ref{sigma_1}) shows that $(\Gamma + I)x=0$. Second,
(\ref{sigma_2}) implies that $\nu_{ij}^Tx=0$ for every $\{i,
j\}\notin E$ and that $\nu_{ij}^Tx + x_i+x_j=0$ for every $\{i,
j\}\in E$. It immediately follows from this last equation that
$\sum_{ (i,j)\in C}\nu_{ij}^Tx=0$ for every cycle $C$, showing
that $x\in \nu(G)^{\perp}$. This ends the proof. \finpr

\begin{cor} Let $G$ be a graph with $n\times n$ adjacency matrix
$\Gamma$. Then \be \log_2|\Sigma_e(G)|\leq n - \mbox{
rank}_{\mathbb{F}_2}(\Gamma +I).\ee
\end{cor}
\begin{cor}
Let $G$ be a graph with adjacency matrix $\Gamma$. If the girth of
$G$ is $\geq 5$ then $|\Sigma_e(G)|=1$.
\end{cor}
{\it Proof: }It was proven in \cite{alg_bouchet} that
$\nu(G)^{\perp}$ is trivial for every graph with girth $\geq
5$.\finpr

\subsection{Invariants of edge-local complementation}

In this section we present some graph invariants under edge-local
equivalence, where we will systematically use theorem 2 to prove
that a given function is an invariant. Of course, there exist many
invariants under edge-local complementation (e.g. all invariants
under local complementation), and we will restrict ourselves to a
handful of interesting ones.

Let $G=(V, E)$ be a graph with $n\times n$ adjacency matrix
$\Gamma$. First, by definition $|L_e(G)|$ is invariant under
edge-local complementation. Second, we show that the space
$\Sigma_e(G)$ is invariant. To see this, let $x\in \Sigma_e(G)$.
Further, let $G'$ be an arbitrary graph edge-locally equivalent to
$G$, and let $\Gamma'$ be the adjacency matrix of $G'$. Then there
exists an $n\times n$ diagonal matrix $A$ such that
$H^A(\Gamma)=\Gamma'$. We then have \be H^{\hat
x}(\Gamma')=H^{\hat x}(H^A(\Gamma))= H^{A+\hat x}(\Gamma) =
H^A(H^{\hat x}(\Gamma))= H^A(\Gamma)=\Gamma',\ee proving that
$x\in\Sigma_e(\Gamma')$. We have proven:
\begin{prop}
Let $G$ be a graph with $n\times n$ adjacency matrix $\Gamma$.
Then the  space $\Sigma_e(G)$ is  invariant under edge-local
complementation.
\end{prop}
\begin{cor}
$|\Sigma_e(G)|$ is invariant under edge-local complementation.
\end{cor}
\begin{cor}
$|\Delta_e(G)|$ is invariant under edge-local complementation.
\end{cor}
{\it Proof:} this follows from the fact that $|L_e(G)|$ and
$|\Sigma_e(G)|$ are invariants. \finpr
\begin{cor}
If $G$ is a graph with adjacency matrix $\Gamma$ satisfying
$\Gamma^2=I$, i.e., $\Gamma$ is an orthogonal matrix, then every
graph edge-locally equivalent to $G$ also has an orthogonal
adjacency matrix.
\end{cor}
{\it Proof: } the result follows by noting that $\Gamma^2=I$ if
and only if $(\Gamma+I)^2=0$, which is equivalent to stating that
the all-ones vector belongs to $\Sigma_e(G)$. The result then
follows from proposition 10. \finpr

 Next we  show that the kernel of $\Gamma + I$ is
invariant under edge-local complementation. To see this, let $G'$
be a graph which is edge-locally equivalent to $G$ and let
$\Gamma'$ be its adjacency matrix, where $\Gamma' = H^A(\Gamma)$
for some diagonal matrix $A$. Let $x\in\mathbb{F}_2^n$ belong to
the kernel of $\Gamma+I$, which is equivalent to stating that $(x,
x)\in V_{\Gamma}$, where the space $V_{\Gamma}$ is defined as in
(\ref{V_Gamma}). It follows that \be (x, x)=H^A(x, x)\in
H^AV_{\Gamma} = V_{\Gamma'},\ee showing that $x$ belongs to the
kernel of $\Gamma'+I$. We have proven:
\begin{prop}
Let $G$ be a graph with $n\times n$ adjacency matrix $\Gamma$.
Then the kernel of the matrix $\Gamma+I$ is  invariant under
edge-local complementation.
\end{prop}
\begin{cor}
The rank of $\Gamma+I$ is  invariant under edge-local
complementation.
\end{cor}

Two vertices $i, j\in V$ of a graph $G=(V, E)$ are called
\emph{twins} if $\{i, j\}\in E$ and $N(i)\setminus\{i\} =
N(j)\setminus\{j\}$, i.e., if these vertices have common
neighborhoods; we can now state the following result.
\begin{cor}
Edge-locally equivalent graphs have the same twins.
\end{cor}
{\it Proof:} let $G$, $G'$ be graphs with adjacency matrices
$\Gamma$, $\Gamma'$ and let $i, j\in V$ be a pair of twins of $G$;
this is equivalent to stating that the $i$th and $j$th column of
$\Gamma+I$ are equal, showing that the vector $x=e_i+e_j$ belongs
to the kernel of this matrix, where $e_i$ ($e_j$) is the $i$th
($j$th) canonical basis vector of $\mathbb{F}_2^n$; it follows
from proposition 11 that $x$ also belongs to the kernel of
$\Gamma'+I$, showing that $i$ and $j$ are also twins of $G'$.
\finpr

\subsection{Interlace polynomial}

For every $k\times k$ matrix $X$ over $\mathbb{F}_2$, we let
$s(X)$ denote the corank of $X$, i.e., the dimension of its
kernel. Letting $G$ be a graph with adjacency matrix $\Gamma$, the
interlace polynomial $q$ of $G$ has the following expansion:
\cite{AH04, interlace_prop}\be q(G, x) = \sum_{\omega\subseteq\{1,
\dots, n\}} (x-1)^{s(\Gamma[\omega])},\ee where by definition
$s(\Gamma[\emptyset])=0$. It is well known \cite{AH04} that
edge-local complementation leaves $q(G, x)$ invariant, such that
edge-locally equivalent graphs have the same interlace polynomial
(in fact, $q$ can be defined through a recursive relation
involving edge-local complementation). We will use the results
presented so far in the present work to gain some insight in the
coefficients of $q(G, x)$. First, it is immediately clear from
(\ref{aantal}) and (\ref{nonsing}) that \be q(G, 1) =
|L_e(G)||\Sigma_e(G)|,\ee linking this evaluation of $q$ with the
number of graphs that are edge-locally equivalent to $G$. In
particular, this result shows that $|\Sigma_e(G)|$ divides $q(G,
1)$. The next proposition shows that a similar result holds for
all coefficients of $q$.

\begin{prop}
Let $G$ be a graph with adjacency matrix $\Gamma$ and let $q(G, x)
= \sum_{i=0}^n a_i x^i$. Then $|\Sigma_e(G)|$ divides every
coefficient $a_i$. In particular, if $\Sigma_e(G)\neq \{0\}$ then
all coefficients $a_i$ are even.
\end{prop}
{\it Proof: }   Let $S\subseteq \mathbb{F}_2^{2n}$ be an arbitrary
$n-$dimensional linear space and let a basis of $S$ be given by
the columns of the matrix, $[Z^T X^T]^T$ 
where $Z$ and $X$ are $n\times n$ matrices. We will say that the
rank of $S$ is equal to $k$ if the rank of $X$ is $k$ (note that
this definition is independent of the chosen basis). Defining
$V_{\Gamma} = \{(\Gamma u, u)\ |\ u\in\mathbb{F}_2^n\}$ as before,
we denote \be L_e^{(k)}(G) =\{ HV_{\Gamma}\ |\ H\in {\cal H}_n,
\mbox{ rank}(HV_{\Gamma})=k\}.\ee Defining \be\Delta_e^{(k)}(G) =
\{H\in{\cal H}_n\ |\ \mbox{ rank}(HV_{\Gamma})=k \},\ee one has
\be |L_e^{(k)}(G)| = \frac{|\Delta_e^{(k)}(G)|}{|\Sigma_e(G)|}.\ee
Moreover, it is easy to verify that $H^A\in \Delta^{(k)}(G)$ if
and only if $s(\Gamma[\omega]) = n-k$, where $A$ is an $n\times n$
diagonal matrix with $\omega=\mbox{ supp}(A)$. Letting \be q(G,
x)= \sum_{k=0}^n b_k (x-1)^k\ee we therefore have \be
|L_e^{(k)}(G)||\Sigma_e(G)|&=&|\Delta_e^{(k)}(G)|\nonumber\\&=&
|\{\omega\subseteq\{1, \dots, n\}\ |\
s(\Gamma[\omega])=n-k\}|\nonumber\\&=&b_{n-k},\ee proving that
$|\Sigma_e(G)|$ divides all coefficients $b_k$. Since \be a_i =
\sum_{k=i}^n (-1)^{i+k}\left( \begin{array}{c} k\\i
\end{array}\right) b_k,\ee it follows that $|\Sigma_e(G)|$ also divides
the coefficients $a_i$.\finpr

\noindent Note that it can efficiently be tested if
$\Sigma_e(G)\neq\{0\}$ for a given graph $G$, since one simply has
to verify whether the matrix (\ref{tilde_Gamma'}) is rank
deficient. This condition is quite strong, such that a typical
graph will not satisfy $\Sigma_e(G)\neq\{0\}$. However, there are
interesting classes of graphs which do satisfy this property. In
the next proposition some sufficient conditions for
$\Sigma_e(G)\neq\{0\}$ to hold are given; this result will be used
below to show that the interlace polynomials of a subclass of
strongly regular graphs have even coefficients.
\begin{prop}
Let $G=(V, E)$ be a graph such that one of the following
situations (i)-(ii) occurs. Then $\Sigma_e(G)\neq\{0\}$.
\begin{itemize}
\item[(i)] $G$ has a pair of twins; \item[(ii)] every vertex of
$G$ has odd degree and $|N(i)\cap N(j)|$ is even for every two
different $i, j\in V$.
\end{itemize}
\end{prop}
{\it Proof: } let $i, j$ be a pair of twins of $G$. Letting
$\tilde\Gamma_i$ and $\tilde\Gamma_j$ denote the $i$th and $j$th
column of $\Gamma+I$ as before, property (i) implies that
$\tilde\Gamma_i=\tilde\Gamma_j$, such that a fortiori
$\tilde\Gamma_i\otimes\tilde\Gamma_i=\tilde\Gamma_j\otimes\tilde\Gamma_j$;
this shows that the matrix (\ref{tilde_Gamma}) is rank deficient,
such that $\Sigma_e(G)$ is nontrivial. If property (ii) holds,
then one has $\Gamma^2=I$ or, equivalently, $(\Gamma+I)^2=0$,
proving that the all-ones vector belongs to $\Sigma_e(G)$. This
yields the result. \finpr

\begin{cor}
Let $G$ be a strongly regular graph with parameters $(n, k, a, c)$
such that $k$ is odd and $a$ and $c$ are even. Then the
coefficients of $q(G, x)$ are all even.
\end{cor}
{\it Proof: } by definition, the degree of each vertex is $k$,
$|N(i)\cap N(j)| = a$ for every $\{i, j\}\in E$, and $|N(i)\cap
N(j)| = c$ for every $\{i, j\}\notin E$. Using proposition 13(ii)
then yields the result.\finpr

\begin{ex} The Clebsch graph $G_{C}$ is the unique strongly
regular graph with parameters $(16, 5, 0, 2)$. Hence, the
interlace polynomials of graphs $G_{C}$ has even coefficients (A
computer calculation showed that $|\Sigma_e(G_C)|=2$).
\end{ex}

\subsection{Inverting an adjacency matrix}

Let $G$ be a graph with adjacency matrix $\Gamma\in GL(n,
\mathbb{F}_2)$, i.e., the matrix $\Gamma$ is nonsingular over
$\mathbb{F}_2$. This is equivalent to stating that the all-ones
vector $d$ belongs to $\Delta(G)$. Moreover, one has \be H^{\hat
d}(\Gamma)=H^I(\Gamma) = (0\cdot\Gamma + I)(I\cdot\Gamma +
0)^{-1}=\Gamma^{-1}. \ee Thus, inverting a nonsingular adjacency
matrix can be realized as an LFT. Consequently, theorem 2 shows
that $\Gamma$ and $\Gamma^{-1}$ are edge-locally equivalent.
Finally, given $\Gamma$ one can calculate $\Gamma^{-1}$ by
applying the following sequence of edge-local complementations:
first, using the method presented in the proof of lemma 1, one
obtains $i_{\alpha}, j_{\alpha}\in \{1 ,\dots, n\}$, where
$\alpha=1, \dots, N$, such that \be f_{i_Nj_N}\dots
f_{i_1j_1}(\Gamma)=I\ee Denoting $e_{\alpha} = \{i_{\alpha},
j_{\alpha}\}$, it follows from proposition 9 that \be
((\Gamma*e_1)*\dots)*e_N = H^I(\Gamma) = \Gamma^{-1},\ee yielding
the desired sequence of edge-local complementation transforming
$\Gamma$ into $\Gamma^{-1}$.

\subsection{Graph states and local Hadamard transformations}

Let $G$ be a graph with $n\times n$ adjacency matrix $\Gamma$ and
let $k_{\Gamma}$ be its associated quadratic form over
$\mathbb{F}_2$, i.e., $k_{\Gamma}(x) =
\sum_{i<j}\Gamma_{ij}x_ix_j$, for every $x=(x_1, \dots,
x_n)\in\mathbb{F}_2^n$. Further, define the canonical basis
vectors $u_0=(1, 0)$ and $u_1 =(0,1)$ in the \emph{real} vector
space $\mathbb{R}^2$. Then one defines the following
$2^n-$dimensional real vector to be the \emph{graph state}
$\psi_G$ associate by $G$: \cite{entgraphstate}
\be \mathbf{\psi}_G
= \frac{1}{2^{n/2}} \sum_{x\in\mathbb{F}_2^n}
(-1)^{k_{\Gamma}(x)}\ u_{x_1}\otimes\dots\otimes u_{x_n}.\ee
In this section we will see that transforming a graph $G$ by
successively applying edge-local complementation or, equivalently,
an LFT $H\in{\cal H}_n$, has a direct translation in terms of a
linear action on the graph state $\mathbf{\psi}_G$. To state this
result, we need some additional definitions. Let
$\omega\subseteq\{1, \dots, n\}$ and define the $2\times2$ real
matrices $\mathbf{H}^{\omega}_i$ by \be \mathbf{H}^{\omega}_i =
\left\{
\begin{array}{cl} H:=\frac{1}{\sqrt{2}}\left[\begin{array}{cc} 1&1\\1&-1\end{array} \right]&\mbox{ if }
i\in\omega\\ I& \mbox{ otherwise},
 \end{array}\right.\ee
where $I$ is here the $2\times2$ real identity matrix; the matrix
$H$ is a $2\times 2$ Hadamard matrix. Furthermore, the $2^n\times
2^n$ real matrix $\mathbf{H}^{\omega}$ is defined by
\be\mathbf{H}^{\omega} = \mathbf{H}^{\omega}_1\otimes\dots\otimes
\mathbf{H}^{\omega}_n.\ee  We call the matrix
$\mathbf{H}^{\omega}$ a local Hadamard matrix\footnote{The term
\emph{local} refers to the fact that $\mathbf{H}^{\omega}$ is an
operator acting on $\mathbb{R}^2\otimes\dots\otimes\mathbb{R}^2$
that can be written as a tensor product of $n$ $2\times 2$
matrices.}. One can now formulate the following result:
\begin{thm}
Let $G$, $G'$ be graphs with $n\times n$ adjacency matrices
$\Gamma$, $\Gamma'$, respectively, and let $A$ be a diagonal
matrix over $\mathbb{F}_2$. Then $H^A(\Gamma)=\Gamma'$ if and only
if \be\mathbf{H}^{supp(A)}\mathbf{\psi}_G\sim
\mathbf{\psi}_{G'},\ee where $\sim$ denotes equality up to a
multiplicative constant.
\end{thm}
This theorem is proven by using the so-called quantum stabilizer
formalism \cite{Gott}. The interested reader is referred to
\cite{doct} for the necessary material to prove the above theorem.
An immediate corollary is the following.
\begin{cor}
Two graphs $G$ and $G'$ are edge-locally equivalent if and only if
there exists a set $\omega\subseteq\{1, \dots, n\}$ such that
$\mathbf{H}^{\omega}\mathbf{\psi}_G \sim \mathbf{\psi}_{G'}$.
\end{cor}
Thus, this result connects edge-local equivalence of graphs with
the action of local Hadamard transformations on graph states.

\section{Conclusion}

In this paper we have studied edge-local equivalence of graphs.
Our main result was an algebraic description of edge-local
equivalence of graphs in terms of linear fractional
transformations (LFTs) over GF(2) of the corresponding adjacency
matrices. Using this equivalent description, we obtained, first, a
polynomial time algorithm to recognize edge-local equivalence of
two arbitrary graphs; the complexity of this algorithm is ${\cal
O}(n^4)$, where $n$ is the number of vertices of the considered
graphs. Second, the description of edge-local equivalence in terms
of LFTs allows to obtain a formula to count the number of graphs
in a class of edge-local equivalence. Third, we considered a
number of invariants of edge-local equivalence. Fourth, we proved
a result concerning the coefficients of the interlace polynomial
of a graph, where we showed that these coefficient are all even
for a class of graphs; this class contains, among others, all
graphs having an orthogonal adjacency matrix (over GF(2)), all
strongly regular graphs with parameters $(n, k, a, c)$, where $k$
is odd and $a$ and $c$ are even, and all graphs having a pair of
twins. Fifth, we showed that a nonsingular adjacency matrix can be
inverted by performing a sequence of edge-local complementations
on the corresponding graph; finally, we considered the connection
between edge-local equivalence of graphs and the action of local
Hadamard matrices on the corresponding graph states.

\appendix

\section{Appendix: Isotropic systems and LFTs}

In this appendix it is shown that the result in theorem 1 is in
fact also present in the work of Bouchet regarding local
complementation and isotropic systems \cite{Bou_iso,
Bou_graph_iso}. In fact, the following analysis will show that the
descriptions of local equivalence of graphs in terms of isotropic
systems and in terms of LFTs are completely equivalent (while they
constitute very different approaches to study the present
subject).

We start be defining isotropic systems \cite{Bou_iso}. Let $K=\{0,
x, y, z\}$ be a two-dimensional vector space over $\mathbb{F}_2$.
There exists a unique inner product $\langle\cdot, \cdot\rangle$
on $K$ satisfying \be \langle a, b \rangle = \left\{
\begin{array}{cl} 1 & \mbox{ if } 0\neq a\neq b\neq 0\nonumber\\ 0
& \mbox{otherwise,} \end{array}\right.\ee for every $a, b\in K$.
Further, let $V$ be a finite set with $n:=|V|$ and consider the
$2n-$dimensional vector space $K^V$ over $\mathbb{F}_2$ consisting
of all $4^n$ functions \be v: i\in V\mapsto v(i)\in K.\ee One
equips the space $K^V$ with the inner product $\langle\cdot,
\cdot\rangle_V$, defined by \be \langle v, w\rangle_V = \sum_{i\in
V} \ \langle v(i), w(i)\rangle,\ee for every $v, w\in K^V$. A
subspace $L\subseteq K^V$ is called \emph{totally isotropic} if
$\langle v, w\rangle=0$ for every $v, w\in L$. An \emph{isotropic
system} is then a pair $S=(L, V)$, where $V$ is a finite set of
cardinality $n$ and $L$ is an $n-$dimensional totally isotropic
subspace of $K^V$.





Let $G =(V, E)$ be a graph on the vertex set $V$. Second, let $a,
b \in K^V$ be two \emph{supplementary vectors}, i.e., $0\neq
a(i)\neq b(i) \neq 0$ for every $i\in V$. Third, for every $v\in
V$ and $\omega\subseteq V$, define $v[\omega]\in K^V$ by \be
v[\omega](i) = \left\{ \begin{array}{cl} v(i) & \mbox{ if } i\in
\omega \\  0 & \mbox{ otherwise. }\end{array}\right.\ee It was
showed in \cite{Bou_graph_iso} that the subspace $L$ of $K^V$
spanned by the vectors \be\label{gen_set}\{a[N(i)] + b[\{i\}]\ |\
i\in V\}\ee is totally isotropic and has dimension $n$, and
therefore $S= (L, V)$ is an isotropic system. The ordered triple
$(G, a, b)$ is called a graphic presentation of the isotropic
system $S$, and $G$ is called a fundamental graph of the system
$S$. Moreover, the following results hold \cite{Bou_graph_iso}:

\begin{thm}
Every isotropic system has a graphic presentation.
\end{thm}

\begin{thm}
Two graphs are fundamental graphs of the same isotropic system if
and only if they are locally equivalent.
\end{thm}

We will now make the connection of the above results with the
theory of LFTs.

First, for every $i\in V$, define the vectors $z^{i},\ x^i\in K^V$
by $z^{i}(i)=z$, $x^i(i)=x$ and $z^{i}(j)= x^i(j) = 0$ for every
$j\in V\setminus\{i\}$. Note that the set $\{z^i, x^i\}_{i\in V}$
is a basis of $K^V$, i.e., every vector $v\in K^V$ can be written
as a linear combination \be v = \sum_{i\in V}\ (v_z)_i z^i +
(v_x)_i x^i,\ee where the coefficients $(v_z)_i,\
(v_x)_i\in\mathbb{F}_2$ are uniquely defined, for every $i\in V$.
Defining $v_z, v_x \in \mathbb{F}_2^V$ by $v_z(i)=(v_z)_i$ and
$v_x(i)=(v_z)_i$, for every $i\in V$, it follows that the mapping
\be\label{iso}\phi: v \in K^V \mapsto \phi(v) = (v_z,
v_x)\in\mathbb{F}_2^{V}\times \mathbb{F}_2^{V}\ee is an
isomorphism between the vector spaces $K^V$ and
$\mathbb{F}_2^{V}\times \mathbb{F}_2^{V}$. We will also identify
the space $\mathbb{F}_2^{V}\times \mathbb{F}_2^{V}$ with the
isomorphic space $\mathbb{F}_2^{2n}$. It is now straightforward to
verify that $\langle v, w \rangle_V = \langle\phi(v),
\phi(w)\rangle$ for every $v, w\in K^V$, where $\langle\cdot,
\cdot\rangle$ is the symplectic inner product on $\mathbb{F}_2^n$
defined in the proof of proposition 1. This leads to the following
result.

\begin{thm} Let $V$ be a finite set with cardinality $n$ and let $\phi$ be the isomorphism defined in (\ref{iso}).
Then $S=(L, V)$ is an isotropic system if and only if $\phi(L)$ is
an $n$-dimensional linear subspace of $\mathbb{F}_2^{2n}$ which is
self-dual w.r.t. the symplectic inner product.
\end{thm}

\noindent Further, let $(G, a, b)$ be a graphic presentation of an
isotropic system $S=(L, V)$. Denote $\phi(a)=(a_z, a_x)$,
$\phi(b)=(b_z, b_x)$ as in (\ref{iso}) and let $\Gamma$ be the
$n\times n$ adjacency matrix of $G$. Finally, let $\Gamma_k$
denote the $k$th column of $\Gamma$ and let $e_k$ denote the $k$th
canonical basis vector of $\mathbb{F}_2^n$, for every $k\in\{1,
\dots, n\}$. With these notations, it is readily verified that the
image of the set (\ref{gen_set}) under the mapping $\phi$ is equal
to the set \be \label{gen_set2} \left\{
\left[ \begin{array}{c} \hat{a_z}\\
\hat{a_x}\end{array}\right] \Gamma_k +
\left[ \begin{array}{c} \hat{b_z}\\
\hat{b_x}\end{array}\right] e_k\ |\ k\in \{1, \dots,
n\}\right\}.\ee  Straightforward manipulations show that the
elements of the set (\ref{gen_set2}) are the columns of the matrix
\be\label{matrix}  \left[
\begin{array}{cc} \hat{a_z} & \hat{b_z}\\
\hat{a_x} & \hat{b_x}\end{array}\right] \left[
\begin{array}{c} \Gamma\\ I\end{array}\right].\ee Recall that the
space $L\subseteq K^V$ is spanned by the set (\ref{gen_set}), as
$(G, a, b)$ is a graphic presentation of $S$. It follows that the
space $\phi(L)\subseteq\mathbb{F}_2^{2n}$ is spanned by the
columns  of the matrix (\ref{matrix}), i.e., \be\phi(L) = \left\{
\left[
\begin{array}{cc} \hat{a_z} & \hat{b_z}\\
\hat{a_x} & \hat{b_x}\end{array}\right] \left[
\begin{array}{c} \Gamma u\\ u\end{array}\right]\ |\ u\in
\mathbb{F}_2^n \right\}\ee
 This shows, in particular, that $\phi(L)$ is the image of the linear space
\be\label{graph_space} V_{\Gamma}:=\{(\Gamma u, u)\ |\ u\in
\mathbb{F}_2^n\}\subseteq\mathbb{F}_2^{2n}\ee  under the mapping
\be \label{Q}Q^{a, b} := \left[
\begin{array}{cc} \hat{a_z} & \hat{b_z}\\
\hat{a_x} & \hat{b_x}\end{array}\right].\ee Furthermore, we have
the following lemma:

\begin{lem}
Let $a, b\in K^V$ and let the corresponding matrix $Q^{a, b}$ be
defined as in (\ref{Q}). Then $a$ and $b$ are supplementary
vectors if and only if $Q^{a, b}\in {\cal C}_n$.
\end{lem}
{\it Proof: } It is sufficient to prove the lemma for the case
where $|V|=1$; the general case follows immediately. If $|V|=1$
then there exist exactly 6 pairs of supplementary vectors, namely
\be\label{supp} (z, x),\ (z, y),\ (x, y),\ (x, z),\ (y,x),\ (y,
z).\ee Also, the isomorphism $\phi$ in this case maps $z\mapsto
(1,0)$, $x\mapsto (0,1)$, $y\mapsto (1,1)$ and one can verify that
the lemma is correct by explicitly constructing the matrices
$Q^{a, b}$ for all 6 pairs of supplementary vectors in
(\ref{supp}). \finpr

\noindent We now arrive at the following result:

\begin{prop}
Let $S=(L, V)$ be an isotropic system and denote $n:=|V|$. Let
$G=(V, E)$ be a graph with adjacency matrix $\Gamma$. Then $(G, a,
b)$ is a graphic presentation of $S$ if and only if $\phi(L) =
Q^{a, b}V_{\Gamma}$, where $Q^{a, b}\in {\cal C}_n$.
\end{prop}

\noindent This leads to the sought-after connection between
isotropic systems and linear fractional transformations.
\begin{cor}
Let $G$ and $G'$ be two graphs on the same vertex set $V$ and let
$\Gamma$ and $\Gamma'$ be the $n\times n$ adjacency matrices of
these graphs. Then $G$ and $G'$ are fundamental graphs of the same
isotropic system if and only if there exist a matrix $Q\in {\cal
C}_n$ such that $\Gamma\in\Delta(Q)$ and $Q(\Gamma)=\Gamma'$.
\end{cor}
{\it Proof: } Let $S=(L, V)$ be the isotropic system such that $G$
and $G'$ are fundamental graphs of $S$. Proposition 14 yields two
matrices $Q_1, Q_2\in{\cal C}_n$ such that $Q_1V_{\Gamma}=\phi(L)
= Q_2V_{\Gamma'}$ and therefore $Q_2^{-1}Q_1 V_{\Gamma} =
V_{\Gamma'}$. Note that, since ${\cal C}_n$ is a group, we have
$Q:= Q_2^{-1}Q_1\in{\cal C}_n$. Furthermore, the identity $Q
V_{\Gamma} = V_{\Gamma'}$ is equivalent to stating that
$\Gamma\in\Delta(Q)$ and $Q(\Gamma) = \Gamma'$ (see the proof of
proposition 1). This completes the proof. \finpr

As a final corollary, we see that theorem 1 now immediately
follows from corollary 10 and theorem 5.

The above analysis shows that the descriptions of local
equivalence of graphs in terms of (graphic presentations of)
isotropic systems and in terms of LFTs are equivalent.

\section{acknowledgments}

MVDN thanks Jeroen Dehaene and Erik Hostens for for many
interesting discussions. This research is supported by several
funding agencies: Research Council KUL: GOA-Mefisto 666,
GOA-Ambiorics, several PhD/postdoc and fellow grants; Flemish
Government: - FWO: PhD/postdoc grants, projects, G.0240.99
(multilinear algebra), G.0407.02 (support vector machines),
G.0197.02 (power islands), G.0141.03 (Identification and
cryptography), G.0491.03 (control for intensive care glycemia),
G.0120.03 (QIT), G.0452.04 (QC), G.0499.04 (robust SVM), research
communities (ICCoS, ANMMM, MLDM); -   AWI: Bil. Int. Collaboration
Hungary/ Poland; - IWT: PhD Grants, GBOU (McKnow) Belgian Federal
Government: Belgian Federal Science Policy Office: IUAP V-22
(Dynamical Systems and Control: Computation, Identification and
Modelling, 2002-2006), PODO-II (CP/01/40: TMS and Sustainibility);
EU: FP5-Quprodis;  ERNSI; Eureka 2063-IMPACT; Eureka 2419-FliTE;
Contract Research/agreements: ISMC/IPCOS, Data4s, TML, Elia, LMS,
IPCOS, Mastercard; QUIPROCONE; QUPRODIS.

\bibliographystyle{unsrt}

\bibliography{localequivgraph}

\end{document}